\documentclass[a4paper,11pt]{article}
\usepackage{amsmath,amssymb,amsfonts,amsthm}
\usepackage{graphicx,color,mathtools}
\numberwithin{equation}{section} 
\mathtoolsset{showonlyrefs=true} 
\usepackage[hidelinks]{hyperref} 

\renewcommand{\title}[1]{\begin{center}{\Large \bf #1}\end{center}}
\renewcommand{\author}[2]{\begin{center}{\large #1\\}
		{\vspace{0.2cm}\small #2}\end{center}}

\newcommand{\email}[1]{\texttt{#1}}

\newtheorem{thm}{Theorem}[section]
\newtheorem{lem}[thm]{Lemma}

\theoremstyle{definition}
\newtheorem{defi}[thm]{Definition}

\newtheorem{problem}[thm]{Problem}

\newtheorem{example}[thm]{Example}
\theoremstyle{remark}
\newtheorem{rmk}[thm]{Remark}

\newcommand{\Hil}{\mathcal{H}}

\newcommand{\slot}{{~\cdot~}}

\DeclareMathOperator{\SO}{SO}
\DeclareMathOperator{\PSL}{PSL}

\newcommand{\C}{\mathcal{C}}
\newcommand{\cF}{\mathcal{F}}

\newcommand{\cU}{\mathcal{U}}

\newcommand{\K}{\mathcal{K}}

\newcommand{\cP}{\mathcal{P}}

\newcommand{\cK}{\mathcal{K}}
\newcommand{\A}{\mathcal{A}}
\newcommand{\B}{\mathcal{B}}
\newcommand{\cI}{\mathcal{I}}

\newcommand{\N}{\mathcal{N}}

\newcommand{\RR}{\mathbb{R}}

\newcommand{\CC}{\mathbb{C}}

\newcommand{\ZZ}{\mathbb{Z}}
\newcommand{\NN}{\mathbb{N}}

\DeclareMathOperator{\Diff}{Diff}

\DeclareMathOperator{\Vir}{{Vir}}

\DeclareMathOperator{\Ad}{Ad}

\DeclareMathOperator{\id}{id}

\DeclareMathOperator{\Aut}{Aut}

\newcommand{\lqq}{\lq\lq}

\usepackage{xspace}
\DeclareRobustCommand{\eg}{e.g.\@\xspace}

\DeclareRobustCommand{\cf}{cf.\@\xspace}
\DeclareRobustCommand{\Cf}{Cf.\@\xspace}
\DeclareRobustCommand{\ie}{i.e.\@\xspace}

\makeatletter
\DeclareRobustCommand{\etc}{%
    \@ifnextchar{.}%
        {etc}%
        {etc.\@\xspace}%
}
\makeatother

\def\u1net{{\A_\RR}}

\DeclareMathOperator*{\QuOp}{QuOp}
\DeclareMathOperator*{\UCP}{UCP}

\def\III{{I\!I\!I}}

\renewcommand{\slot}{\,\cdot\,}

\date{}

\begin{document}
	
	\begin{flushright}
	\end{flushright}
	
	\title{\Large{\bf{Quantum Operations in Algebraic QFT}}}
	
	\author{{\bf Luca Giorgetti}$^a$ }
	{
		$^a$ Dipartimento di Matematica, Universit\`a di Roma Tor Vergata, \\
		Via della Ricerca Scientifica, 1, I-00133 Roma, Italy \\ 
		
		\email{giorgett@mat.uniroma2.it}} 
	
	Conformal Quantum Field Theories (CFT) in 1 or 1+1 spacetime dimensions (respectively called chiral and full CFTs) admit several
\lqq axiomatic" (mathematically rigorous and model-independent) formulations. In this note,
we deal with the von Neumann algebraic formulation due to Haag and Kastler \cite{HaagBook}, mainly restricted to the chiral CFT setting \cite{CKLW18}.
Irrespectively of the chosen formulation, one can ask the questions: given a theory $\mathcal{A}$, how many and
which are the possible extensions $\mathcal{B} \supset \mathcal{A}$ or subtheories 
$\mathcal{B} \subset \mathcal{A}$? How to construct and classify them, and study their properties? 
Extensions are typically described in the language of \emph{algebra objects} in the braided tensor category
of representations of $\mathcal{A}$, while subtheories require different ideas.

In this paper, we review recent structural results on the study of subtheories in the von Neumann algebraic formulation (conformal subnets) 
of a given chiral CFT (conformal net), \cite{Bis17}, \cite{BDG21}, \cite{BDG22}, \cite{BDG23}. Furthermore, building on \cite{BDG23}, 
we provide a \lqq quantum Galois theory" for conformal nets analogous to the one for Vertex Operator Algebras (VOA) \cite{DoMa97}, \cite{DoMa99}.
We also outline the case of 3+1 dimensional Algebraic Quantum Field Theories (AQFT). 

The aforementioned results make use of families of (extreme) vacuum state preserving unital completely positive maps acting
on the net of von Neumann algebras, hereafter called \emph{quantum operations}. These are natural generalizations of the 
ordinary vacuum preserving gauge automorphisms, hence they play the role of \lqq generalized global gauge symmetries".
Quantum operations suffice to describe \emph{all} possible conformal subnets of a given conformal net with the same central charge.

\section{Conformal nets}

Conformal nets are the von Neumann algebraic description of \emph{chiral}, \ie, 1-dimensional (living on a single lightray of a 1+1-dimensional Minkowski spacetime) CFT. We identify the lightray with $\RR$ and we consider for convenience its one-point compactification $S^1$ (the unit circle). 
Denote by $\cI$ be the set of non-empty non-dense open intervals on $S^{1}$, and by $I' := (S^{1}\setminus I)^{\circ} \in \cI$ the interior of the complement of $I\in\cI$. Let $\mathrm{PSL}(2,\RR) := \mathrm{SL}(2,\RR)/\{\pm 1\}$ be the M\"obius group acting by fractional linear transformations on $S^1$. Let also $\B(\Hil)$ and $\cU(\Hil)$ be respectively the algebra of bounded operators on $\Hil$ and its unitary subgroup.
We refer \eg, to \cite[Chapter 3]{CKLW18} and \cite{LonLec} for self-contained introductions and references.

\begin{defi}
A {\bf M\"obius covariant net} on $S^{1}$ is a triple $(\A, U, \Omega)$ consisting of a family of von Neumann algebras $\A=\left\{\A(I) \subset\B(\Hil): I\in\cI\right\}$ acting on a common separable Hilbert space $\Hil$, a strongly continuous unitary representation $U : \mathrm{PSL}(2,\RR) \to \cU(\Hil)$ and a unit vector $\Omega \in \Hil$, satisfying the following properties:
\begin{enumerate}
\item[(i)] \textbf{Isotony}: $\A(I_{1})\subset\A(I_{2})$, if $I_{1}\subset I_{2}$, $I_{1},I_{2}\in \cI$.
\item[(ii)] \textbf{Locality}: $\A(I_{1})\subset\A(I_{2})'$, if $I_{1}\cap I_{2}=\emptyset$, $I_{1},I_{2}\in \cI$.
\item[(iii)] \textbf{M\"obius covariance}: if $I\in\cI$, $g\in \mathrm{PSL}(2,\RR)$,
\begin{align}
U(g)\A(I)U(g)^{-1}=\A(gI).
\end{align}
\item[(iv)] \textbf{Positivity of energy}: $U$ has positive energy. Namely, the conformal Hamiltonian (the generator of the one-parameter rotation subgroup of $\mathrm{PSL}(2,\RR)$) has non-negative spectrum.
\item[(v)] \textbf{Vacuum vector}: $\Omega$ is the unique (up to a phase) vector with the property
$U(g)\Omega=\Omega$ for every $g\in \mathrm{PSL}(2,\RR)$, and vectors of the form $x\Omega$, $x\in\bigvee_{I\in\cI}\A(I)$, are dense in $\Hil$.
\end{enumerate}

Here $\bigvee_{I\in\cI}\A(I)$ denotes the von Neumann algebra generated on $\Hil$, and $\A(I)'$ denotes the \emph{commutant} of $\A(I)$ in $\B(\Hil)$, namely $\A(I)' := \{x\in\B(\Hil): xy = yx, y\in\A(I)\}$. In particular, by von Neumann's bicommutant theorem, $\A(I) = \A(I)''$. For an introduction, see, \eg, \cite{JonLec}. The $\A(I)$ will be referred to as the \emph{local algebras}, and $\Hil$ as the \emph{vacuum Hilbert space} of $\A$.
\end{defi}

The M\"obius group plays the same role of the Poincar\'e group for QFTs in 3+1 dimensions. In the low-dimensional conformal setting, the quantum fields typically enjoy a much bigger (infinite dimensional) symmetry group: the group of orientation preserving diffeomorphisms of the unit circle, $\Diff_+(S^1)$.

\begin{defi}\label{def:subnet}
A M\"obius covariant net $(\A,U,\Omega)$ will be called a \textbf{conformal net} on $S^1$, if it satisfies in addition:
\begin{enumerate}
\item[(vi)] $U$ extends to a strongly continuous projective unitary representation of $\Diff_+(S^1)$, again denoted by $U$, such that for every $I\in\cI$,
\begin{align}
U(\gamma)\A(I)U(\gamma)^{-1}&=\A(\gamma I), \quad\gamma\in\Diff_+(S^1),\\
U(\gamma)xU(\gamma)^{-1}&=x,\quad x\in \A(I), \gamma\in\Diff_+(I^\prime),
\end{align}
where $\Diff_+(I^\prime)$ denotes the subgroup of orientation preserving diffeomorphisms of $S^1$ that are \emph{localized} in $I'$, namely $\gamma\in\Diff_+(S^1)$ such that $\gamma(z)=z$ for all $z\in I$.
\end{enumerate}
\end{defi}

Conformal nets allow for a model-independent treatment of chiral CFT, and many examples have been constructed and classified in this setting, see, \eg, \cite{KaLo04}. Mathematically, the local algebras turn out to be type $\III_1$ factors, the most non-commutative infinite-dimensional type of von Neumann algebras (with no non-zero faithful trace/tracial weight) \cite{ConnesBook} and with trivial center (equal to $\CC 1$, the scalar operators on $\Hil$). The vacuum vector $\Omega$ is cyclic and separating for each local algebra $\A(I)$ (the so-called \emph{Reeh--Schlieder property}). Namely, $\A(I)\Omega$ is a dense subspace of $\Hil$, and $x\Omega = 0$, $x\in\A(I)$, implies $x = 0$.

\begin{rmk}
The terminology \lqq factor" comes from the easiest (finite-dimensional) examples of such von Neumann algebras: the full $n\times n$ complex matrix algebras $M_n(\CC)$. If $M_n(\CC)$ acts (faithfully) on some $\CC^{m}$, $m\in\NN$, then the space \lqq factorizes" as $\CC^m \cong \CC^n \otimes \CC^k$ for some multiplicity index $k\in\NN$, where $M_n(\CC)$ acts on the first $\CC^n$ and $(M_n(\CC))' \cong M_k(\CC)$ (the commutant is taken in $\CC^m$) acts on the second $\CC^k$. Not all factors are like this \cite{MuvN36}, namely $\N = \N''\subset\B(\Hil)$ with trivial center $\N\cap\N' = \CC1$, \ie, $\N\vee\N' = \B(\Hil)$, does not imply a tensor \lqq factorization" of the Hilbert space $\Hil \cong \Hil_1 \bar\otimes \Hil_2$, such that $\N$ acts faithfully on $\Hil_1$ and $\N'$ acts faithfully on $\Hil_2$.
\end{rmk}

\section{Conformal subnets}\label{sec:confsubnets}

Once a chiral CFT is given, \eg, in the form of a conformal net $\A$, the problem of studying and classifying its \emph{extensions} can be phrased in terms of algebra-like objects in the (braided tensor) representation category of $\A$, \cite{LoRe95}. More precisely, \emph{relatively local} extensions $\B \supset \A$, \ie, such that $\B(I_1) \subset \A(I_2)'$ if $I_1 \cap I_2 = \emptyset$, $I_{1},I_{2}\in \cI$, with \textbf{finite Jones index} \cite{Jon83}, \cite{Kos86} (\ie, \lqq finite relative size") are described by certain C$^*$-Frobenius algebra objects (also called Q-systems \cite{Lon94}, \cf \cite{BKLR15}, \cite{Gio22}) in the representation category of $\A$. The \emph{locality} of $\B$ itself can be characterized in terms of the associated Q-system as a suitable \lqq commutativity constraint".
In the Vertex Operator Algebra (VOA) formulation of chiral and full (1+1-dimensional) CFT, see, \eg, respectively \cite{KacBook}, \cite{HuKo07}, the same structure of commutative algebras in the module category of a given theory turns out to describe its \lqq finite order" local extensions, \cite{KiOs02}, \cite{HKL15}, \cite{Kon07}.
In the finite index/finite order and chiral CFT case, the two theories of extensions for conformal nets and for unitary VOAs \cite{CKLW18} have been unified in \cite{Gui22}, \cite{CGGH22}.
In the conformal net setting, it has also been generalized to the \textbf{infinite index discrete} case, \cite{DVGi18}, \cf \cite{ILP98}, \cite{Car04}, and some seemingly new models of full CFT have recently been constructed in \cite{AGT23}. While for the general \textbf{semidiscrete} case \cite{DVGi18}, \cf \cite{FiIs99}, it is not clear how to proceed. The study of \lqq infinite order" extensions has recently been tackled in the VOA framework as well, \cite{CMY22}.

The study of \emph{subtheories} of a given theory $\A$ is a completely different problem, no matter if finite or infinite index. 
The structure of subtheories is not anymore governed by the representation theory of $\A$ (\eg, the so-called \emph{holomorphic} models have no non-vacuum irreducible representations but may have many inequivalent subtheories $\B\subset\A$).

\begin{defi}\label{confsubnetdef}
A \textbf{conformal subnet} of a conformal net $(\A,U,\Omega)$ is a family $\B= \{\B(I) : I\in\cI\}$ of non-trivial von Neumann algebras acting on the same space $\Hil$ of $\A$ such that:
\begin{enumerate}
\item[(i)] $\B(I) \subset \A(I)$ for every $I \in \cI$.
\item[(ii)]  $U(g)\B(I)U(g)^{-1} = \B(g I)\,$ for every $I \in \cI, g \in \PSL(2,\RR)$.
\item[(iii)] $\B(I_1)\subset \B(I_2)$ for every $I_1\subset I_2$, $I_1,I_2\in\cI$.
\end{enumerate}
\end{defi}

\begin{rmk}
A conformal subnet $\B\subset\A$ fulfills also diffeomorphism covariance with respect to the same $U$, see, \eg, \cite[Section 3.4]{CKLW18}, \cite[Section 2]{BDG23}:
\begin{align}
U(\gamma)\B(I)U(\gamma)^{-1} &= \B(\gamma I), \quad I\in\cI, \gamma\in\Diff_+(S^1),\\
U(\gamma)xU(\gamma)^{-1} &= x,\quad x\in \B(I), \gamma\in\Diff_+(I^\prime).
\end{align}
\end{rmk}

The second equation above follows from the second equation in Definition \ref{def:subnet}, just because $\B(I) \subset \A(I)$.
From this second equation, by the so-called \emph{Haag duality property} of the net $\A$ on $\Hil$ (\ie, $\A(I)' = \A(I')$ for $I\in\cI$, which is a consequence of the cyclic and separating property of the vacuum vector $\Omega$ for $\A$ on $\Hil$, \cite{FrJoe96}), it follows that $U(\gamma) \in \A(I)$ if $\gamma\in\Diff_+(I)$. Note however that Haag duality for $\B$ does not hold on $\Hil$, hence $U(\gamma) \in \B(I')'$ if $\gamma\in\Diff_+(I)$, but $U(\gamma)$ may or may not belong to $\B(I)$. Let
\begin{align}
\Vir_\A(I) := \{U(\gamma): \gamma\in\Diff_+(I)\}^{\prime\prime}, \quad I \in \cI,
\end{align}
and $(\Vir_\A, U, \Omega)$ be the \emph{Virasoro subnet} of $(\A,U,\Omega)$.
By the previous discussion, $\Vir_\A$ may or may not be a subnet of $\B$.
The Virasoro nets are completely determined by a parameter $c$ called \emph{central charge} (whose possible values are quantized below the value $1$, namely $c = 0$, $c = 1 - \frac{6}{m(m+1)}$ for $m\geq 3$, $m\in\NN$, and continuous above $1$), \cite{BPZ84}, \cite{FQS85}.

The conformal inclusion $\Vir_\A \subset \A$ is automatically \emph{irreducible}, \cite{Car04} (\ie, $\Vir_\A(I)' \cap \A(I) = \CC 1$ for every $I\in\cI$), and it is \emph{minimal}, \cite{Car98}, in the sense that it has no non-trivial conformal subnet.
A conformal inclusion $\B \subset \A$ is \emph{irreducible} (\ie, $\B(I)' \cap \A(I) = \CC 1$ for every $I\in\cI$) if and only if the condition mentioned above $U(\gamma) \in \B(I)$, if $\gamma\in\Diff_+(I)$, is fulfilled. In this case,
$$\Vir_\A \subset \B \subset\A,$$
moreover, $\Vir_\A = \Vir_\B$, and $\A$ and $\B$ have the same central charge. As in \cite{BDG23}, this is the type of subnets we deal with. They are called \textbf{conformal inclusions} in the physics literature. Moreover, all finite index inclusions are necessarily irreducible, \cite{Lon03}, hence conformal. \Cf \cite{CGH19} for a family of examples of conformal subnets with infinite index that are not conformal inclusions in the sense specified above.

The most studied families of examples of subtheories come from group fixed points subnets $\B = \A^G \subset \A$, where $G$ acts on $\A$ by \emph{automorphisms}, the so-called \emph{group orbifolds} \cite{DLM96}, \cite{DoMa97}, \cite{Xu00}, \cite{Xu05}, \cite{DRX17}. In the conformal net setting, automorphisms are defined as follows, see, \eg, \cite[Section 3.3 and 5.3]{CKLW18}:

\begin{defi}\label{def:Aut}
Let $(\A,U,\Omega)$ be a conformal net. An \textbf{automorphism} of $\A$ is a map $\alpha = \Ad V (\equiv V \slot V^{-1})$ for some unitary $V\in\cU(\Hil)$ (\ie, $V^*V = 1 = VV^*$) such that:
\begin{enumerate}
\item[(i)] \textbf{Compatibility}: $V \A(I) V^{-1} = \A(I)$ for every $I\in\cI$.
\item[(ii)] \textbf{Vacuum preserving}: $V \Omega = \Omega$.
\end{enumerate}
Denote by $\Aut(\A)$ the set of all automorphisms of $\A$, which clearly forms a group under composition and inversion of unital *-automorphisms of each local algebra. 
\end{defi}

\begin{rmk}
The following \emph{covariance condition} on $\alpha$ (or $V$) is a consequence of (i) and (ii), namely (iii): $U(\gamma) V U(\gamma)^{-1} = V$ for every $\gamma\in\Diff_+(S^1)$, see, \eg, \cite{GaFr93}.
\end{rmk}

\begin{rmk}
$\Aut(\A)$ turns out to be \emph{compact} in the strong operator topology on the $V$, as a combination of results in \cite{DoLo84} and of the \emph{split property}, \cite{MTW18}. 
\end{rmk}

\begin{rmk}\label{rmk:Gorbifold}
For every closed subgroup $G \subset \Aut(\A)$, the corresponding orbifold subnet $\A^G \subset \A$ gives rise to a conformal inclusion, see \cite{Car99}, \cite{Xu05}. For group orbifolds the index is either finite (and $[\A:\A^G] = |G|$, the order of $G$) or infinite (if $G$ is infinite). In the latter case, the inclusion $\A^G \subset \A$ is necessarily discrete \cite{Car04}, \cite{DVGi18}.
\end{rmk}

We recall below some examples that are in a sense orthogonal to group orbifolds:

\begin{example}
The (local) conformal nets $\A$ with central charge $c<1$ are classified in \cite{KaLo04}.  
There are the trivial inclusions $\Vir_\A = \A$ with index $[\A:\Vir_\A] = 1$ and any central charge $c = 1 - \frac{6}{m(m+1)}$ for any $m\geq 2$, $m\in\NN$. Moreover, there are index $[\A:\Vir_\A] = 2$ inclusions, almost trivial in the sense that they are $\ZZ_2$-fixed points subtheories $\Vir_\A = \A^{\ZZ_2} \subset \A$ (or, the other way around, $\ZZ_2$-simple current extensions), either with central charge given by $m = 4n+1$ or by $m = 4n+2$, for any $n \in \NN$. Lastly, there are four \lqq exceptional cases", two with index $3 + \sqrt{3}$ and respectively $m=11$, $m=12$, the other two with index approximately $19.48$ and respectively $m=29$, $m=30$. Let us take one of the first two exceptional cases, namely $\Vir_\A \subset \A$ with $m=11$, hence $c=\frac{21}{22}$, and index $[\A: \Vir_\A] = 3+\sqrt{3}$. Clearly, there cannot be a finite group $G\subset\Aut(\A)$ (of order $3 + \sqrt{3}$) such that $\Vir_\A = \A^G \subset \A$. Moreover, $\Aut(\A)$ turns out to be trivial for this conformal net $\A$ (a fact that can be deduced by the very classification of \cite{KaLo04}, or by subfactor theoretical arguments). There are also less extreme examples where both inclusions $\B \subset \A^{\Aut(\A)} \subset \A$ are non-trivial.
\end{example}

We shall later see that arbitrary irreducible conformal subnets (no matter if finite or infinite index and discrete or not) can be described as the fixed points under \emph{quantum operations} (see Section \ref{sec:QuOp}), \ie, as \emph{generalized orbifolds} in the terminology of \cite{Bis17}.

\section{Quantum operations}\label{sec:QuOp}

In order to describe arbitrary subtheories in a model-independent way, in the operator algebraic setting, we introduced in \cite{BDG23} the notion of \emph{quantum operation} on a conformal net. The main novelty is that unital completely positive (UCP) maps of each local algebra into itself are used instead of automorphisms, see, \eg, \cite{PauBook}, \cite{Tak1}. We proved in \cite{BDG23} that quantum operations suffice to describe all subtheories of a given chiral CFT. In this section, we review some of the results in \cite{BDG23}.

\begin{defi}
Let $(\A,U,\Omega)$ be a conformal net. A \textbf{unital completely positive map} on $\A$ is a collection of normal faithful unital completely positive maps:
\begin{align}
\phi=\{\phi_I : \A(I)\rightarrow\A(I),\quad I\in\cI\},
\end{align}
such that:
\begin{enumerate}
\item[(i)] \textbf{Compatibility}: $\phi_{I_2} \restriction_{\A(I_1)}=\phi_{I_1}$ for every $I_1\subset I_2$, $I_1,I_2\in \cI$.
\item[(ii)] \textbf{Vacuum preserving}: $\omega_I = \omega_I \circ \phi_I$ on $\A(I)$, where $\omega_I := (\Omega, \slot \Omega)$ is the vacuum state restricted to $\A(I)$, for every $I\in\cI$.
\item[(iii)] \textbf{M\"obius covariance}: $\Ad U(g) \circ \phi_I\circ \Ad U(g)^{-1} = \phi_{g I}$ for every $I\in\cI$, $g\in\mathrm{PSL}(2,\RR)$.
\item[(iv)] \textbf{$\Vir$-fixing}: $\phi_I(U(\gamma)) = U(\gamma)$ for every $I\in\cI$, $\gamma \in \Diff_+(I)$.
\end{enumerate}
Denote by $\UCP(\A)$ the set of all unital completely positive maps on $\A$, which is a convex set as convex combinations preserve the UCP property and the conditions (i)-(iv) above.
\end{defi}

\begin{rmk}
Condition (iv), namely $\phi_I \restriction_{\Vir_\A(I)}\, = \id$, is equivalent to $\Vir_\A(I)$-bimodularity of $\phi_I$, namely $\phi_I(axb) = a \phi_I(x) b$ if $x\in\A(I)$ and $a,b\in\Vir_\A(I)$.
\end{rmk}

\begin{rmk}
Assuming in addition that the inclusion $\Vir_\A \subset \A$ is finite index (\eg, finite group of automorphisms fixed points) or discrete (\eg, compact group of automorphisms fixed points, as described in the previous section), then both covariance conditions (iii) and (iv) are a consequence of (i) and (ii), as shown in \cite[Section 6]{BDG23}.
\emph{Question:} Do the conditions (iii) and (iv) follow from (i) and (ii) in full generality, \ie, only having semidiscreteness of the inclusion at our disposal?
\end{rmk}

The group of automorphisms $\Aut(\A)$ (Definition \ref{def:Aut}) is contained in $\UCP(\A)$. Indeed, $\alpha \in \Aut(\A)$ restricted to each $\A(I)$ defines an element $\phi \in \UCP(\A)$ by $\phi_I := \alpha_{\restriction\A(I)}$, as the unital *-automorphisms are special examples of UCP maps.
Vice versa, if some $\phi\in\UCP(\A)$ happens to be invertible in $\UCP(\A)$ (or just multiplicative, \ie, $\phi_I(xy) = \phi_I(x) \phi_I(y)$ if $x,y\in\A(I)$), then necessarily $\phi\in\Aut(\A)$.

\begin{rmk}
Similarly to the implementing unitary $V$ of a vacuum preserving automorphism $\alpha$, one can consider the \emph{standard implementation} of $\phi\in\UCP(\A)$, see \cite[Section 4]{BDG23}. Namely, let $V_\phi x \Omega := \phi_I(x) \Omega$ for every $x\in\A(I)$, for a fixed $I\in\cI$, where note that the vectors $x \Omega$ are dense in $\Hil$ by the Reeh--Schlieder property. It turns out that $V_{\phi}$ is bounded with norm 1 (by the Kadison--Schwarz inequality, \cite{Kad52}), hence it can be defined on the whole $\Hil$. It is easy to show that $V_\phi$ does not depend on $I$ by (i). Moreover, $U(\gamma) V_\phi U(\gamma)^{-1} = V_\phi$ for every $\gamma\in\Diff_+(S^1)$ by (iii). Clearly, $V_\phi \Omega = \Omega$.

One can show that the standard implementation $V_\phi$ is unitary (or isometric, \ie, $V_\phi^* V_\phi = 1$) if and only if $\phi$ is in $\Aut(\A)$ (or multiplicative, hence in $\Aut(\A)$ as already mentioned).\footnote{For every $\phi\in\UCP(\A)$ (or in $\QuOp(\A)$), the Hilbert space adjoint $V_\phi^*$ corresponds to another element $\phi^\sharp\in\UCP(\A)$ (or in $\QuOp(\A)$), defined as the \emph{$\Omega$-adjoint UCP map} of $\phi_I$ on each $\A(I)$. See \cite{BDG23}.}
\end{rmk}

\begin{defi}
Denote also by $\QuOp(\A)$ the subset of \emph{extreme points} of $\UCP(\A)$ (in the sense of convex sets). We call \textbf{quantum operations} on $\A$ the elements of $\QuOp(\A)$.
\end{defi}

\begin{rmk}
One reason for considering the extreme points of $\UCP(\A)$, instead of the whole $\UCP(\A)$, is that $\Aut(\A)$ sits into $\QuOp(\A)$ as well, \cite{BDG23}. Another reason is that, in the special case of group orbifolds: if $\Vir_\A = \A^G \subset \A$, then $\QuOp(\A) = \Aut(\A) = G$ and $\UCP(\A) \cong P(G)$, the convex set of probability Radon measures on $G$, as a combination of results in \cite{BDG21} and \cite{BDG23}. Note that $G$ is identified with the extreme points of $P(G)$ by considering the normalized Dirac delta measures centered at the elements of $G$.
\end{rmk}

The convex set $\UCP(\A)$ turns out to be \emph{compact} in a natural locally convex vector space (Hausdorff) topology, induced by the pointwise ultraweak operator topology on each local algebra, see \cite[Theorem 4.10]{BDG23}. From this, together with the existence of a vacuum preserving \emph{conditional expectation} from $\A$ \emph{onto} $\Vir_\A$ in $\UCP(\A)$ (\cf the discussion in \cite[Remark 2.5]{BDG23}), by applying the Krein--Milman theorem we get:

\begin{thm}\label{thm:virisgenorbi}
Let $(\A,U,\Omega)$ be a conformal net. Then $\Vir_\A = \A^{\QuOp(\A)}$. 

In particular, every irreducible conformal subnet $\B\subset\A$ is the fixed points under a subset of quantum operations $\B = \A^S$, for some $S \subset \QuOp(\A)$.
\end{thm}

\begin{rmk}\label{rmk:QuOpiscpthypifdiscrete}
On the one hand, $\Aut(\A)$ is always a group and compact in a natural topology, as we already said. However, as we have seen in Section \ref{sec:confsubnets}, it does not suffice to describe arbitrary subnets of $\A$, not even in the finite index case.
On the other hand, a priori we know nothing about the compactness of $\QuOp(\A)$ itself, as the extreme points need \emph{not} be closed, nor Borel, in general. Moreover, we do not know what is the algebraic structure of $\QuOp(\A)$, except when $\QuOp(\A) = \Aut(\A)$, as the composition of extreme UCP maps need \emph{not} be extreme. 

If $\Vir_\A \subset \A$ happens to be discrete or finite index, it turns out that $\QuOp(\A)$ is a \emph{compact hypergroup}, see \cite{BDG21}, \cite{BDG23}, and references therein.\footnote{A \emph{hypergroup} is an abstract convolution space generalizing the notion of group (in a classical, not quantum, way). The operation is defined on probability measures, instead of on points, and the convolution of two extreme probability measures (Dirac delta measures) need not be extreme anymore. It is always extreme precisely when the hypergroup is a group, see, \eg, \cite{Jew75}.}
\end{rmk}

In \cite{BDG23}, we also define and study \textbf{relative quantum operations} on $\A$ with respect to an irreducible conformal subnet $\B \subset \A$, hence with $\Vir_\B = \Vir_\A$, but possibly $\B \neq \Vir_\A$. This is because $\Vir_\A \subset \A$ is basically never discrete, except when it has finite index, see, \eg, \cite{Fre94}, \cite{Reh94}, \cite{Car03}. Also, relative quantum operations allow one to study \emph{intermediate} subnets.

\begin{defi}
Let $\QuOp(\A|\B) := \{\phi\in\QuOp(\A) : {\phi_I}\restriction_{\B(I)} = \id, \forall I\in\cI\}$.
Let also $\UCP(\A|\B) := \{\phi\in\UCP(\A) : {\phi_I}\restriction_{\B(I)} = \id, \forall I\in\cI\}$. Then one can show that $\QuOp(\A|\B)$ coincides with the subset of extreme points of $\UCP(\A|\B)$.
\end{defi}

\emph{Question:} Given a conformal subnet $\B\subset\A$, what are the topological and algebraic structures of $\QuOp(\A) (= \QuOp(\A|\Vir_\A))$ and $\QuOp(\A|\B)$, in the general (semidiscrete) case?

\section{Galois correspondence for conformal subnets}

In \cite{DLM96}, \cite{DoMa97}, \cite{DoMa99}, \cite{HMT99}, a \lqq quantum Galois theory" for VOAs is established. Namely, given a simple VOA $V$ and a (finite or compact Lie) group of automorphisms $G \subset \Aut(V)$. Let $V^G$ be the $G$-fixed points subVOA of $V$. Then \cite[Theorem 3]{DoMa99} establishes a one-to-one correspondence between the closed Lie subgroups of $G$ and the orthogonally complemented, \cite[Definition 2]{DoMa99}, simple subVOAs of $V$ that contain $V^G$. Building on \cite{BDG23} and \cite{ILP98}, we give below (Theorem \ref{thm:quantumgalois}) the analogous statement for conformal nets and compact (not necessarily Lie) groups of automorphisms $G \subset \Aut(\A)$.

\begin{rmk}\label{rmk:UCPlocisglob}
If $\B \subset \A$ has finite index or is irreducible and infinite index discrete, the second main technical result of \cite{BDG23}, namely \cite[Theorem 6.4]{BDG23}, states that every $\B(I)$-fixing vacuum preserving UCP map $\A(I) \to \A(I)$ (in particular, every $\B(I)$-fixing vacuum preserving automorphism of $\A(I)$), for a fixed $I\in\cI$, can be extended to a coherent family of UCP maps defining an element of $\UCP(\A|\B)$. Similarly, every $\B(I)$-fixing vacuum preserving automorphism of $\A(I)$ can be extended to an element of $\Aut(\A|\B)$. Moreover, the map $\phi\in\UCP(\A|\B) \mapsto \phi_I$ is an affine homeomorphism.
\end{rmk}

Recall now from Remark \ref{rmk:Gorbifold} that $\A^G \subset \A$ is always irreducible and discrete (with finite or infinite index depending on the size of $G$) whenever $G\subset\Aut(\A)$ is a compact group of automorphism of $\A$. Consequently, as a special case of \cite[Theorem 7.4]{BDG23}, we get:

\begin{thm}\label{thm:quantumgalois}
Let $(\A, U, \Omega)$ be a conformal net and $G\subset\Aut(\A)$ be a compact group of automorphisms of $\A$. Then there is a one-to-one correspondence between the closed subgroups $H \subset G$ and the conformal subnets $\C \subset \A$ that contain $\A^G$, given by $H \mapsto \A^H$.
\end{thm}

The previous Theorem \ref{thm:quantumgalois} also follows from the aforementioned extension result for UCP maps (or automorphisms) \cite[Theorem 6.4]{BDG23} combined with the well-known group-subgroup Galois correspondence for subfactors, see \cite[Theorem 3.15]{ILP98}.

\begin{rmk}
Whenever $\B\subset \A$ is irreducible and discrete or finite index, in \cite[Theorem 7.4, 7.6]{BDG23} a generalization of Theorem \ref{thm:quantumgalois} is given by means of closed sub\emph{hyper}groups of relative quantum operations $\QuOp(\A|\B)$. The latter turns out to have itself the structure of a compact hypergroup (\cf the last comment in Remark \ref{rmk:QuOpiscpthypifdiscrete} where $\B = \Vir_\A$), it acts on $\A$ by definition and it gives back the subnet $\B = \A^{\QuOp(\A|\B)}$ by a relative version of Theorem \ref{thm:virisgenorbi}. 
\emph{Question:} Given a conformal net $(\A, U, \Omega)$ and given a compact hypergroup $K \subset \QuOp(\A)$ (thus acting by vacuum preserving UCP maps, \cf \cite[Section 5]{BDG21}), is the fixed point inclusion $\A^K \subset \A$ automatically irreducible and discrete, as for compact groups of automorphisms $G\subset \Aut(\A) ( \subset \QuOp(\A) )$?
\end{rmk}

\section{Quantum operations on 3+1-dimensional AQFTs}

The operator algebraic formulation of QFT was originally intended for Poincar\'e covariant theories in 3+1-dimensional Minkowski spacetime. See, \eg, \cite{HaKa64}, \cite{HaagBook}, \cite{BauBook}. In this more physical situation, the theory of finite or infinite index extensions from \cite{LoRe95}, \cite{DVGi18} can still be applied. However, in the 3+1-dimensional situation, much more can be said after the work of Doplicher and Roberts \cite{DoRo90}, based on \cite{DoRo89}, \cite{DoRo89-1} and \cite{DHR71}, \cite{DHR74}. Namely, to every net of local observables $\A = \{\A(O)\}$, indexed by the open non-empty bounded regions $O \subset \RR^{3+1}$ (describing the local measurements of a QFT) one can associate a graded-local \emph{canonical field net} $\cF = \{\cF(O)\}$ that extends $\A$, and a canonical compact group of automorphisms $G\subset\Aut(\cF)$ (describing the \lqq global gauge group symmetries" of $\cF$) such that $\A = \cF^G$. Using this result (which has no analogue in 1 or 1+1 dimensions, \cite{Mue03}), the intermediate local nets $\A \subset \B \subset \cF$ \cite{CDR01} and the Poincar\'e covariant subnets $\B \subset \A$ have been classified \cite{CaCo01}, \cite{CaCo05}, see in particular \cite[Theorem 5.2]{CaCo05} and \cf \cite{CaCo01-proc}. Below, we define quantum operations on 3+1 dimensional nets of local observables and we review the previously mentioned classification results in this language.

Let $\cK$ be the set of open non-empty double cones in Minkowski spacetime $\RR^{3+1}$ and let $\cP := \SO(3,1)^{\uparrow}_+$ be the  proper orthochronous Poincar\'e group.

\begin{defi}
A {\bf Poincar\'e covariant net of local observables} (or just \textbf{local net}) on Minkowski spacetime $\RR^{3+1}$ is a triple $(\A, U, \Omega)$ consisting of a family of von Neumann algebras $\A=\left\{\A(O) \subset\B(\Hil): O\in\cK\right\}$ acting on a common separable Hilbert space $\Hil$, a strongly continuous unitary representation $U : \cP \to \cU(\Hil)$ and a unit vector $\Omega \in \Hil$, satisfying the following properties:
\begin{enumerate}
\item[(i)] \textbf{Isotony}: $\A(O_{1})\subset\A(O_{2})$, if $O_{1}\subset O_{2}$, $O_{1},O_{2}\in \cK$.
\item[(ii)] \textbf{Locality}: $\A(O_{1})\subset\A(O_{2})'$, if $O_{1},O_{2}\in \cK$ are space-like separated in $\RR^{3+1}$.
\item[(iii)] \textbf{Poincar\'e covariance}: if $O\in\cK$, $g\in \cP$,
\begin{align}
U(g)\A(O)U(g)^{-1}=\A(gO).
\end{align}
\item[(iv)] \textbf{Positivity of energy}: $U$ has positive energy. Namely, the joint spectrum of the generators of spacetime translations (a subgroup of $\cP$ isomorphic to $\RR^4$) is contained in the closure of the open forward light cone.
\item[(v)] \textbf{Vacuum vector}: $\Omega$ is the unique (up to a phase) vector with the property
$U(g)\Omega=\Omega$ for every spacetime translation $g$ in $\cP$.
\end{enumerate}

Further additional properties are not automatic in the 3+1-dimensional situation, unlike in the chiral conformal setting, and they are often included in the definition of local net. For example, the Reeh--Schlieder property: $\Omega$ is cyclic (and separating) for each $\A(O)$, $O\in\cK$. The Haag duality property: $\A(O)' = \bigvee_{\tilde O \in \cK, \tilde O \subset O'} \A(\tilde O)$, where $O'$ is the interior of the space-like complement of $O$ in $\RR^{3+1}$. Further assumptions (namely: CPT covariance, the Bisognano--Wichmann property, and the split property) are also usually included. See, \eg, \cite{HaagBook}, \cite{CaCo01}, \cite{CaCo01-proc}, and references therein.
\end{defi}

\begin{defi}
A \textbf{Poincar\'e covariant subnet} of $(\A,U,\Omega)$ is a family $\B= \{\B(O) : O\in\cK\}$ of non-trivial von Neumann algebras acting on the same space $\Hil$ of $\A$ such that:
\begin{enumerate}
\item[(i)] $\B(O) \subset \A(O)$ for every $O \in \cK$.
\item[(ii)]  $U(g)\B(O)U(g)^{-1} = \B(g O)\,$ for every $O \in \cK, g \in \cP$.
\item[(iii)] $\B(O_1)\subset \B(O_2)$ for every $O_1\subset O_2$, $O_1,O_2\in\cK$.
\end{enumerate}
\end{defi}

As for $\A$, we assume that $\B$ fulfills Haag duality. We also assume that the inclusion $\B \subset \A$ is \emph{irreducible}, or better that $\B$ is \emph{full} in $\A$, in the terminology of \cite[Section 2]{CaCo01}. As in the chiral CFT case \cite{BDG23}, we define quantum operations for 3+1-dimensional QFTs:

\begin{defi}
Let $(\A,U,\Omega)$ be a Poincar\'e covariant net of local observables. A \textbf{unital completely positive map} on $\A$ is a collection of normal faithful unital completely positive maps:
\begin{align}
\phi=\{\phi_O : \A(O)\rightarrow\A(O),\quad O\in\cK\},
\end{align}
such that:
\begin{enumerate}
\item[(i)] \textbf{Compatibility}: $\phi_{O_2} \restriction_{\A(O_1)}=\phi_{O_1}$ for every $O_1\subset O_2$, $O_1,O_2\in \cK$.
\item[(ii)] \textbf{Vacuum preserving}: $\omega_O = \omega_O \circ \phi_O$ on $\A(O)$, where $\omega_O := (\Omega, \slot \Omega)$ is the vacuum state restricted to $\A(O)$, for every $O\in\cK$.
\item[(iii)] \textbf{Poincar\'e covariance}: $\Ad U(g) \circ \phi_O \circ \Ad U(g)^{-1} = \phi_{g O}$ for every $O\in\cK$, $g\in\cP$.
\end{enumerate}
Denote by $\UCP(\A)$ the set of all unital completely positive maps on $\A$, and by $\QuOp(\A)$ its extreme points, which we call \textbf{quantum operations} on $\A$.
\end{defi}

\begin{rmk}
In the 3+1-dimensional context, there is no obvious canonical analogue of the Virasoro net sitting irreducibly and minimally inside every local net, \cf \cite{BDL86}, \cite{Con95}, \cite{CaCo01-proc} where an irreducible local subnet of $\A$ generated by the implementation of the translations via $U$ is considered. \emph{Question:} Does however the covariance condition (iii) follow from (i) and (ii)? Here note that the \emph{discreteness} assumption on an inclusion of local nets $\B\subset\A$ is regarded as natural in the 3+1-dimensional context, \cf \cite{LoRe95}, \cite{CDR01}, \cite{CaCo01}, \cite{CaCo05}, as the irreducible DHR representations of a local net typically have finite statistical/C$^*$-tensor categorical dimension \cite{DHR71}, \cite{LoRo97}. 
\end{rmk}

With the analogous definition of vacuum preserving automorphisms of $\A$ in the 3+1-dimensional context, see, \eg, \cite[Section 3]{DoRo90}, \cite{HaagBook}, and with the same argument used \cite{BDG23} for chiral CFTs, namely using the proof of \cite[Lemma 4.49]{BDG21}, it follows that:

\begin{lem}
Let $(\A,U,\Omega)$ be a Poincar\'e covariant net of local observables. Then $\Aut(\A) \subset \QuOp(\A)$.
\end{lem}

One can presumably prove the analogous results as for the chiral conformal net case. For instance: prove that $\UCP(\A)$ is compact in a natural topology induced by the pointwise ultraweak operator convergence on the UCP maps $\phi_O$. More importantly, in the 3+1-dimensional situation, by the classification of subsystems due to \cite{CDR01}, \cite{CaCo01}, \cite{CaCo05}, in particular \cite[Theorem 5.2]{CaCo05}, we propose the following:

\begin{problem}
Let $(\A,U,\Omega)$ be a Poincar\'e covariant net of local observables, and $\B\subset\A$ a (Haag dual, full) Poincar\'e covariant subnet. Then $\QuOp(\A|\B) \cong G /\!\!/ H$ (a \emph{double coset hypergroup}\footnote{If $H$ is a normal subgroup of $G$, then $G /\!\!/ H$ is a group (the quotient group), otherwise it is a prototypical example of compact hypergroup (a double quotient hypergroup), see, \eg, \cite[Chapter 1.5]{BlHeBook}.}) for a compact group $G$ (obtained as the Doplicher--Roberts global gauge group of the subtheory $\B$) and a closed subgroup $H \subset G$.
\end{problem}

\begin{rmk}
In particular, in the 3+1-dimensional context, we conjecture that $\QuOp(\A)$ is always compact and a hypergroup (as close as possible to a group, in some sense).
\end{rmk}

For a fixed region $O\in\K$, namely, for a single subfactor $\B(O) \subset \A(O)$ (as opposed to the whole net of subfactors $\B\subset\A$), a positive solution of the previous problem follows from the results in \cite{BDG21}, see, in particular, \cite[Section 9.2 and Example 9.26]{BDG21}. One way of solving the problem in general would be to show the analogue of \cite[Theorem 6.4]{BDG23} for Poincar\'e covariant local nets.

\subsubsection*{Acknowledgements}
I thank Maria Stella Adamo, Sebastiano Carpi, Simone Del Vecchio, Roberto Longo, Alessio Ranallo, Yoh Tanimoto for stimulating and helpful discussions. I also thank Roberto Longo for drawing my attention to \cite{Fre94} and Ajit Iqbal Singh for email correspondence and for drawing my attention to \cite{Sin96}, where families and hypergroups of UCP maps on operator algebras are studied in greater generality.

This work was supported by the European Union's Horizon 2020 research and innovation programme H2020-MSCA-IF-2017 under Grant Agreement n.\! 795151 \emph{Beyond Rationality in Algebraic CFT: mathematical structures and models}, acronym: beyondRCFT. It has also been supported by the \emph{MIUR Excellence Department Project MatMod@TOV}, CUP E83C23000330006, awarded to the Department of Mathematics, University of Rome Tor Vergata, by the University of Rome Tor Vergata funding OAQM, CUP E83C22001800005, and by GNAMPA--INdAM.

Invited contribution to the conference proceedings of \emph{Functional Analysis, Approximation Theory and Numerical Analysis}, FAATNA20$>$22,
Matera, Italy, July 5--8, 2022.

	
	
\end{document}